\documentclass[journal]{IEEEtran}

\usepackage{myStyleIEEE}
\usepackage[margin=0.5in]{geometry}
\allowdisplaybreaks

\IEEEoverridecommandlockouts

\begin{document}

\title{Improving Dynamic Performance of Low-Inertia Systems through Eigensensitivity Optimization}

\renewcommand{\theenumi}{\alph{enumi}}

\newcommand{\uros}[1]{\textcolor{magenta}{UM: #1}}
\newcommand{\vaggelis}[1]{\textcolor{blue}{EV #1}}
\newcommand{\petros}[1]{\textcolor{red}{PA #1}}
\newcommand{\ashwin}[1]{\textcolor{pPurple}{AV: #1}}

\author{Ashwin~Venkatraman,~\IEEEmembership{Student~Member,~IEEE,}
        Uros~Markovic,~\IEEEmembership{Member,~IEEE,}
        Dmitry~Shchetinin,~\IEEEmembership{Member,~IEEE,}
        Evangelos~Vrettos,~\IEEEmembership{Member,~IEEE,}
        Petros~Aristidou,~\IEEEmembership{Senior~Member,~IEEE,}
        Gabriela~Hug,~\IEEEmembership{Senior~Member,~IEEE}\vspace{-0.5cm}% <-this % stops a space\vspace{-0.35cm}
}

\maketitle
\IEEEpeerreviewmaketitle

\begin{abstract}
An increasing penetration of renewable generation has led to reduced levels of rotational inertia and damping in the system. The consequences are higher vulnerability to disturbances and deterioration of the dynamic response of the system. To overcome these challenges, novel converter control schemes that provide virtual inertia and damping have been introduced, which raises the question of optimal distribution of such devices throughout the network. This paper presents a framework for performance-based allocation of virtual inertia and damping to the converter-interfaced generators in a low-inertia system. This is achieved through an iterative, eigensensitivity-based optimization algorithm that determines the optimal controller gains. Two conceptually different problem formulations are presented and validated on a 3-area, 12-bus test system.
\end{abstract}

%INDEX TERMS

\section{Introduction} \label{sec:intro}

\lettrine[lines=2]{T}{he electric} power system is currently undergoing a major transition integrating large shares of distributed generation interfaced via Power Electronic (PE) converters. This is accompanied by the phase-out of conventional Synchronous Generators (SGs) leading to a loss of rotational inertia. Such developments have serious effects on system dynamics, especially in terms of performance \cite{tielens2016relevance}, such as larger frequency excursions and Rate-of-Change-of-Frequency (RoCoF) after disturbances \cite{EirGrid2012,ERCOTwind}. To address the underlying issues, different technologies and control techniques have been proposed, primarily in terms of restoring the lost inertia by emulating a Synchronous Machine (SM)-like behavior through converter control. A number of control schemes have been designed to make power converters behave as closely as possible to SMs \cite{Beck2007,Zhong2014,AJD18}. These schemes can vary in the level of detail and complexity, but they all rely on replicating the dynamic behavior of SGs, therefore providing virtual inertia and damping. Moreover, all such strategies require energy on the DC side acting as converter's equivalent to the missing kinetic energy of the rotor; either in the form of energy storage (e.g., batteries, flywheels, supercapacitors), or by employing a power source with available kinetic energy (e.g., wind turbines and diesel generators) \cite{Van2010,Torres2009}. 

Synthetic inertia and damping are thus becoming design parameters of the power system and could potentially serve as a foundation for ancillary service \cite{EirGrid2012,ERCOTwind}. Therefore, a natural next step is to quantify the required amount of virtual inertia and damping for a particular system. However, it is not clear yet based on which performance metrics such quantification should take place. Traditionally, the total aggregated inertia and damping in the system were used as the main metrics used for measuring system resilience~\cite{NordicReport2016}. However, the authors of \cite{BKP-SB-FD:17,Ulbig2014} show that, in addition to aggregate inertia and damping, the distribution of these parameters in the system could be of crucial importance, with spatially heterogeneous inertia profiles resulting in worse dynamic response after a disturbance. Other performance metrics such as frequency nadir, RoCoF and minimum damping ratio are also commonly used in the literature~\cite{Rakhshani2016,Borsche2015}. In contrast, the problem of optimal tuning and placement of Virtual Synchronous Machine (VSM) control gains pertaining to virtual inertia and damping has recently been tackled from the perspective of system norms \cite{BKP-SB-FD:17,xiong2018stability,ademola2018optimal,Poolla2019,Pirani2017,Paganini2017,Mesanovic16}, namely using the $\mathcal{H}_2$ norm in \cite{BKP-SB-FD:17,ademola2018optimal,GrossIREP,xiong2018stability,Poolla2019}, $\mathcal{H}_\infty$ norm in \cite{Pirani2017}, $\mathcal{L}_2$ and $\mathcal{L}_\infty$ norms in \cite{Paganini2017}, and all of the above approaches in \cite{Mesanovic16}.

The existing literature suggests two main directions to approach this problem. On the one hand, \cite{Borsche2015} analyzes the sensitivity of eigenvalues with respect to inertia and damping, and subsequently maximizes the critical damping ratio of the system while ensuring that frequency overshoot is limited. Although such formulation is non-convex by nature, it can be linearized and solved numerically in an iterative fashion. Nevertheless, the study considers an oversimplified representation of the system, especially in terms of modeling of converter-interfaced generation, and employs a sequence of approximations to obtain an estimate of the frequency metrics of interest. On the other hand, \cite{BKP-SB-FD:17,xiong2018stability,ademola2018optimal,Poolla2019} use the $\mathcal{H}_2$ norm as a measure of network coherency and characterization of the system frequency response, as well as for quantifying the VSM control effort. In particular, \cite{BKP-SB-FD:17} and \cite{xiong2018stability} aim at improving the frequency response of the system by finding the inertia distribution that minimizes the $\mathcal{H}_2$ norm, with the second study specifically focusing on a network with high penetration of wind farms employing doubly-fed induction generators. However, neither of the studies include the damping constant as a controllable parameter. This drawback is resolved in \cite{ademola2018optimal,Poolla2019}, where the authors argue that the performance metrics such as the damping ratio and RoCoF are not sufficient for quantifying the robustness of the system. They instead employ the $\mathcal{H}_2$ norm not only to characterize the system response, but also to quantify the required control effort. Furthermore, \cite{ademola2018optimal} incorporates simplified Virtual Inertia (VI) devices operating in the grid-following mode as a feedback control loop, whereas \cite{Poolla2019} also includes the grid-forming VI implementation. While providing certain theoretical guarantees, the use of $\mathcal{H}_2$ norm is limited by the underlying assumption of an impulse disturbance, which is not necessarily the case in power systems (e.g., load change and generation outage yield a step-like change in active power). Another challenge is the computational burden to solve the underlying Lyapunov equations for a detailed low-inertia system. This was not an obstacle in \cite{BKP-SB-FD:17,xiong2018stability,ademola2018optimal,Poolla2019} due to the simplistic representation of system dynamics based on the swing equation. However, such models are not sufficient to characterize the dynamic interactions present in a realistic system with high inverter penetration \cite{UrosStability}, thus raising concerns in terms of applicability of system norms.

The work proposed in this paper aims to combine the techniques and insightful conclusions from the existing literature and provides a methodology for placement of virtual inertia and damping in a high-fidelity low-inertia system. To this end, a detailed dynamic model of such system presented in \cite{UrosStability} is used for optimal VSM control design (i.e., control tuning) and assessing the system performance. Similarly to~\cite{Rakhshani2016,Borsche2015}, we employ an iterative, eigensensitivity-based optimization framework to determine the optimal incremental allocation of virtual gains at each step. However, unlike the approximations made in \cite{Borsche2015}, we include exact analytical expression for relevant frequency metrics previously derived in \cite{UrosLQR}. Moreover, improvements in terms of the multi-objective nature of the problem, computational efficiency and adaptive step-size adjustments are also presented. Additionally, motivated by \cite{ademola2018optimal,Poolla2019}, we include different VI implementations as well as the $\mathcal{H}_2$ and $\mathcal{H}_\infty$ system norms in the analysis. While computationally intractable within the iterative algorithm, the two norms are studied and taken into consideration when evaluating the system performance. As a result, we obtain a computationally inexpensive and scalable problem formulation that takes into account all relevant aspects of the dynamic system response. To the best of the knowledge of the authors, such comprehensive and multifaceted approach has not been proposed in the literature thus far.

The rest of the paper is organized as follows. Section~\ref{sec:problem} discusses the principles of eigensensitivity optimization and the different performance metrics and problem objectives. The mathematical formulation of the problem is presented in Section~\ref{sec:formulation}, followed by case studies in Section~\ref{sec:res}. Finally, Section~\ref{sec:conclusion} discusses the future work and concludes the paper. 

\section{Theoretical Preliminaries} \label{sec:problem}

\subsection{Dynamic Properties of Low-Inertia Systems} \label{subsec:dyn_prop}
The level of inertia and damping present in the network largely influences the small-signal stability and frequency dynamics after a disturbance. %With the increased penetration of converter-interfaced generation, the amount of inertia and damping reduces and the system becomes increasingly vulnerable to disturbances such as load fluctuations and generation outages. 
A few commonly used metrics for assessing the dynamic behavior of the system are listed and briefly discussed below:
\begin{itemize}
    \item \textit{Small-signal stability}: Defined as the ability of the system to maintain synchronism when subject to a small disturbance. It was previously shown in \cite{UrosStability} that the system dominated by both grid-forming and grid-following power converters can face small-signal instability under insufficient levels of virtual inertia and damping (i.e., small VSM control gains).
    \item \textit{Damping ratio}: Describes how fast the oscillations in the system die out. A higher damping ratio increases system resilience.  
    \item \textit{Frequency nadir}: Represents the maximum deviation of frequency from a nominal value after a disturbance. Frequency nadir is a nonlinear function of both inertia and damping, as will be shown in Section~\ref{subsec:freq_metrics}.
    \item \textit{Rate-of-Change-of-Frequency}: Describes the maximum rate at which the system frequency changes and usually corresponds to the instantaneous RoCoF value after a disturbance. Unlike frequency nadir, the RoCoF is solely a function of system inertia.
\end{itemize}

\subsection{Eigensensitivity Optimization Principles} \label{subsec:eig_principle}
Improving the worst-case damping ratio of all modes in the system is important for ensuring an acceptable dynamic response. The damping ratios are functions of, among other parameters, inertia and damping constants of both synchronous and converter-based generators. However, being functions of system eigenvalues, the sensitivities of damping ratios to respective parameters are highly nonlinear and could result in a complex optimization problem \cite{Borsche2015}. This section provides a brief introduction into the computation of such sensitivities and how they can be incorporated into a sequential iterative algorithm for improving the damping ratios.

The general state-space representation of a linearized system is given by
\begin{align}
    \dot{x}&=Ax+Bu, \label{eq:state-space1}\\
    y&=Cx+Du, \label{eq:state-space2}
\end{align}
where $x\in\R^n$, $y\in\R^m$ and $u\in\R^p$ are the respective state, output and control input vectors, and $A\in\R^{n\times n}$, $B\in\R^{n\times p}$, $C\in\R^{m\times n}$ and $D\in\R^{m\times p}$ are the state-space matrices. Let $\lambda_i\in\mathbb{C}$ denote the $i$-th eigenvalue of the system and $\sigma_i\in\R$ and $\omega_i\in\R$ its real and imaginary parts. The right and left eigenvectors $u_i, v_i\in\mathbb{C}^n$ of $\lambda_i$ are then described as
\begin{align}
    A u_i = \lambda_i u_i, \\
    v_i^\mathsf{T} A =v_i^\mathsf{T} \lambda_i,
\end{align}
whereas the damping ratio $\zeta_i\in\R$ is defined as
\begin{align}
    \zeta_i \coloneqq \frac{-\sigma_i}{\sqrt{\sigma_i ^2 + \omega_i ^2}}.
\end{align}
This value is positive for stable modes, zero for oscillatory modes and negative for unstable modes. Moreover, the sensitivity of eigenvalue $\lambda_i$ with respect to a parameter $\alpha\in\R$ can be expressed as a function of eigenvectors, i.e., 
\begin{align}
    \frac{\partial \lambda_i}{\partial \alpha} = \frac{\partial \sigma_i}{\partial \alpha} + j\frac{\partial \omega_i}{\partial \alpha} = v_i^\mathsf{T} \frac{\partial A}{\partial \alpha} u_i,
\end{align}
which in turn can be used to compute the sensitivity of damping ratio $\zeta_i$ with respect to the parameter $\alpha$, as follows:
\begin{align} \label{eq:zeta_sensitivities}
    \frac{\partial \zeta_i}{\partial \alpha} = \frac{\partial}{\partial \alpha}\left(\frac{-\sigma_i}{\sqrt{\sigma_i ^2 + \omega_i ^2}}\right) = \omega_i \frac{\sigma_i \frac{\partial \omega_i}{\partial \alpha} - \omega_i \frac{\partial \sigma_i}{\partial \alpha}}{\left(\sigma_i^2 + \omega_i^2 \right)^\frac{3}{2}}.
\end{align}

Since the underlying sensitivities are nonlinear, the task of maximizing the damping ratio is performed using an iterative approach \cite{Borsche2015}. More precisely, the sensitivities are obtained at the start of the iteration using \eqref{eq:zeta_sensitivities}, from which the new parameter values that maximize the damping ratios are computed. Such iterative procedure yields
\begin{align}
    \zeta^{\nu+1}_i =\zeta^{\nu}_i + {\frac{\partial \zeta^{\nu}_i}{\partial \alpha}} \left( \alpha^{\nu+1} - \alpha^{\nu}\right),
\end{align}
where $\nu\in\N_{0}$ denotes the iteration step and $\alpha^{\nu}$ and $\alpha^{\nu+1}$, i.e., $\zeta^{\nu}$ and $\zeta^{\nu+1}$, represent the old and new values of parameters and damping ratios, respectively. The updated damping ratios and the corresponding sensitivities are subsequently used in the next iteration step, described in more detail in Section~\ref{subsec:soln_strategy}. 

\subsection{Applicability of System Norms} \label{subsec:apply}

Apart from the metrics presented in Section~\ref{subsec:dyn_prop}, system norms such as $\mathcal{H}_2$ and $\mathcal{H}_\infty$ provide a measure of the magnitude of the system output after a disturbance. The system output can include performance outputs such as frequency stability and energy of the control effort, thereby making system norms a useful tool for optimization. In general, the $\mathcal{H}_2$ norm measures the energy of the system's impulse response, whereas the $\mathcal{H}_\infty$ norm represents the peak gain from the disturbance to the output~\cite{Zhou1996}. By defining a suitable performance output, the energy metrics of the VSM control effort can be directly considered in the $\mathcal{H}_2$ framework as the overall output energy~\cite{GrossIREP}. Nevertheless, both system norms have drawbacks when it comes to computation and applicability to a detailed model of a low-inertia system.

The $\mathcal{H}_2$ norm can be computed using the controllability and observability Gramians of the system. While tractable for small-scale systems, the computation of two Lyapunov functions pertaining to the controllability and observability Gramians becomes numerically intensive on larger systems. Furthermore, the $\mathcal{H}_2$ norm quantifies the system performance when subjected to an impulse disturbance, which is not applicable to more common disturbances in the power system such as a loss of generator or a load demand change. On the other hand, the $\mathcal{H}_\infty$ norm is not restrictive in terms of the nature of the disturbance signal. However, unlike the $\mathcal{H}_2$ norm, its computation cannot be expressed in a concise analytical form and requires the use of iterative algorithms such as those presented in \cite{boyd1989bisection}, \cite{Bruinsma1990}. The bisection method in \cite{boyd1989bisection} has similar computational drawbacks as the computation of the $\mathcal{H}_2$ norm. On the other hand, the less computationally intensive method provided in \cite{Bruinsma1990} is not suitable for our purposes, since the eigensensitivity framework presented in Section~\ref{subsec:eig_principle} requires the computation of norm sensitivities to decision variables (specifically inertia and damping). It should be noted that the computational effort for obtaining the sensitivities of aforementioned norms to system parameters is also an additional obstacle for explicitly including them into the eigensensitivity-based problem.

Finally, an important question to consider is the practicality of minimizing the system norms. Indeed, the main concern for system operators is to ensure that the frequency metrics such as RoCoF and frequency nadir are within the limits prescribed by the system operator in order to prevent false triggering of protection and load shedding schemes~\cite{entsoeFreq}. Meeting these requirements is sufficient for providing reliable and safe operation and any further improvement of frequency response (i.e., minimization of speed and magnitude of frequency deviation) is not necessarily of value to the operator. By minimizing system norms significantly tighter bounds are imposed on these frequency metrics, for which the required control effort and cost may not be justified. Moreover, optimizing system norms is not directly correlated with the damping ratio and does not guarantee achieving sufficient damping of oscillatory modes. Even though it is hard to prove any formal relationship between the $\mathcal{H}_\infty$ norm of a linear time-invariant system and the damping ratio of its eigenvalues, the intuition suggests that by improving the latter one could also reduce the former. The results presented in Section~\ref{sec:res} also support these claims. Therefore, the method proposed in this paper prioritizes the improvement of the worst-case damping ratio while ensuring that the frequency constraints are met, which simultaneously leads to a compelling reduction of system norms.

\section{Problem Formulation} \label{sec:formulation}

This section describes the proposed optimization problem. The overview of the dynamic model and analytical expressions for frequency constraints are presented, followed by two different problem formulations, i.e., multi-step and uniform. In particular, the multi-step approach comprises three steps, each of them with a different objective function and set of constraints, whereas the uniform formulation combines all of these objectives and constraints into a single optimization problem. 

\subsection{Dynamic Model of Low-Inertia System} \label{subsec:sys_model}

We consider a state-of-the-art Voltage Source Converter (VSC) control scheme previously described in \cite{UrosISGTeurope}, where the outer control loop consists of VSM-based active and droop-based reactive power controllers providing the output voltage angle and magnitude reference by adjusting the predefined setpoints according to a measured power imbalance. Subsequently, the reference voltage vector signal is passed through a virtual impedance block as well as the inner control loop consisting of cascaded voltage and current controllers. The DC voltage is controlled through a DC current source and a PI controller. To detect the system frequency at the connection terminal, a synchronization unit in the form of a phased-locked loop is included in the model of grid-following VSC units. For more details, the reader is referred to \cite{UrosISGTeurope}. 

For SGs we consider a traditional $6^\text{th}$-order round rotor generator model equipped with a prime mover and a \textit{TGOV1} governor. Furthermore, the Automatic Voltage Regulator (AVR) based on a simplified excitation system \textit{SEXS} is incorporated, together with a \textit{PSS1A} power system stabilizer. Detailed control configuration and tuning parameters are provided in \cite{Kundur1994,entsoeGen,UrosStability}. The SG is interfaced through a transformer to the grid and modeled in the $(dq)$-frame defined by its synchronous speed. Finally, transmission network dynamics are also included in the model, with short transmission lines modeled as $\pi$-sections and long transmission lines represented by distributed line parameters. More details on the mathematical formulation and dynamic performance of the employed model can be found in \cite{UrosStability}.

\subsection{System Frequency Metrics} \label{subsec:freq_metrics}
 
To incorporate frequency constraints into the optimization problem, the first step is to obtain the corresponding analytical expressions in terms of decision variables (i.e., parameters of interest), namely inertia and damping. The detailed derivation of frequency nadir and RoCoF metrics after a step disturbance of magnitude $\Delta P\in\R$ is presented in \cite{UrosLQR}. The underlying expressions in SI are given by
\begin{subequations} \label{eq:freq_terms}
\begin{align}
    \dot f_\mathrm{max} &\coloneqq - f_0 \frac{\Delta P}{M}, \label{eq:RoCoF_exp}\\
    \Delta f_\mathrm{max} &\coloneqq - f_0 \frac{\Delta P}{D+R_g} \Bigg(1+ \sqrt{\frac{T(R_g - F_g)}{M}} e^{-\zeta_s \omega_n t_m}\Bigg). \label{eq:Nadir_exp}
\end{align}
\end{subequations}
Here, $\dot f_\mathrm{max}\in\R$ and $\Delta f_\mathrm{max}\in\R$ are the maximum values of RoCoF and frequency nadir, $f_0=\SI{50}{\hertz}$ represents the nominal frequency, $M\in\R_{\geq0}$ and $D\in\R_{\geq0}$ denote the weighted system averages of inertia and damping, $R_g\in\R_{\geq0}$ is the average inverse droop control gain (effectively corresponding to damping) and $F_g\in\R_{\geq0}$ represents the fraction of the total power generated by the high pressure turbines of the SGs. The time instance of the frequency nadir is given by
\begin{equation}
    t_m \coloneqq \frac{1}{\omega_n \sqrt{1-\zeta_s^2}} \tan^{-1}{\frac{\omega_n \sqrt{1-\zeta_s^2}}{\zeta_s \omega_n - T^{-1}}} \label{eq:tm}
\end{equation}
with $\zeta_s\in\R_{>0}$ and $\omega_n\in\R_{>0}$ representing the damping ratio and the natural frequency of the system response:
\begin{equation}
    \zeta_s \coloneqq \frac{M+T\left(D+F_g\right)}{2 \sqrt{MT\left(D+R_g\right)}}\,, \quad 
    \omega_n \coloneqq \sqrt{\frac{D+R_g}{M+T}}.
\end{equation}
An additional constraint is introduced (see \cite{UrosLQR} for more details) to ensure that the time instance of frequency nadir is positive, i.e., $t_m>0$, which corresponds to $\frac{M}{T} - F_g < D$.

These expressions are highly dependent on the aggregate system inertia and damping, and the limits enforced on frequency metrics in \eqref{eq:freq_terms} can be translated into bounds on $M$ and $D$. However, frequency nadir in \eqref{eq:Nadir_exp} is a nonlinear function of system parameters. To incorporate such constraint into the iterative linear program, with inertia and damping being decision variables, a first-order Taylor approximation is employed:
\begin{align}
    \Delta f_\mathrm{max}^{\nu+1} = \Delta f_\mathrm{max}^{\nu} + \frac{\partial \Delta f_\mathrm{max}^\nu}{\partial M} \Delta M^{\nu +1} + \frac{\partial \Delta f_\mathrm{max}^\nu}{\partial D} \Delta D^{\nu +1},
\end{align}
where $\Delta f_\mathrm{max}^{\nu}$ and $\Delta f_\mathrm{max}^{\nu+1}$ are the values of frequency nadir in iteration $\nu$ and $\nu+1$, and $\Delta M^{\nu +1}\in\R$ and $\Delta D^{\nu +1}\in\R$ are the respective updates of system inertia and damping at each iteration step.

\subsection{Multi-Step Optimization Problem} \label{subsec:multistep}
 
Let $\mathcal{N}\subset\N$ be the set of network buses, $\mathcal{K}\subseteq\mathcal{N}$ and $\mathcal{J}\subseteq\mathcal{N}$ represent the subset of nodes with synchronous and converter-interfaced generation respectively, and $n_k=|\mathcal{K}|, n_j=|\mathcal{J}|$. The per-unit inertia and damping constants of unit $j$ are described by $m_j\in\R_{>0}$ and $d_j\in\R_{>0}$, whereas their incremental changes computed at iteration step $\nu$ are denoted by $\Delta m^{\nu+1}\in\R^{n_j}$ and $\Delta d^{\nu+1}\in\R^{n_j}$. Note that we only consider the virtual inertia and damping gains of converter-interfaced generators as decision variables (i.e., control gains to be tuned), and the parameters of SGs remain intact. Moreover, $\sigma_i^\nu\in\R$ and $\zeta_i^\nu\in\R$ reflect the real part and the damping ratio of the $i^\mathrm{th}$ mode at iteration step $\nu$. 

To address all performance metrics listed in Section~\ref{subsec:dyn_prop}, we propose a sequential procedure comprising three consecutive optimization problems. This multi-step approach first addresses the small-signal stability of the system, followed by improving the worst-case damping ratio. Finally, the inertia and damping are redistributed across the system such that the total amount of additional control effort is minimized. 

\subsubsection{Step 1: Ensuring Small-Signal Stability}
The first step aims at restoring the small-signal stability of the system by ensuring that the real parts of all eigenvalues become negative, i.e., $\sigma_i<0,\forall i\in \N_{\leq n}$. Let us define $\Phi_1 \coloneqq [m^{\nu+1}, d^{\nu+1}, \Delta m^{\nu+1}, \Delta d^{\nu+1}]^\mathsf{T}$ as the vector of optimization variables at each iteration step. The problem can be formulated as
\begin{subequations} \label{opti_stability}
\begin{alignat}{3}
    & \underset{\Phi_1, \sigma^{\nu+1}}{\min} && \sigma_\mathrm{max} \label{1:obj}\\
    & \;\mathrm{s.t.} \quad && \forall j \in \mathcal{J}, \forall i \in \N_{\leq n}, \nonumber \\
    & && \sigma^{\nu+1}_i = \sigma^{\nu}_i + \sum_{j \in \mathcal{J}}{\frac{\partial \sigma^{\nu}_i}{\partial d_j}} \Delta d_j^{\nu+1} + \sum_{j \in \mathcal{J}}{\frac{\partial \sigma^{\nu}_i}{\partial m_j}} \Delta m_j^{\nu+1}, \label{1:sigma_iter}\\
    & && \sigma^{\nu+1}_i \leq \sigma_\mathrm{max}, \label{1:sigma_lim}\\
    & && \ushort{d}_j \leq d_j^{\nu+1} \leq \widebar{d}_j, \label{1:k_lim}\\ 
    & && \ushort{m}_j \leq m_j^{\nu+1} \leq \widebar{m}_j, \label{1:m_lim}\\ 
    & && \Delta \ushort{d}_j\phi^d_j \leq \Delta d_j^{\nu+1} \leq \Delta \widebar{d}_j \phi^d_j, \label{1:step_k}\\
    & && \Delta \ushort{m}_j\phi^m_j \leq \Delta m_j^{\nu+1} \leq \Delta \widebar{m}_j\phi^m_j, \label{1:step_m}\\
    & && \Delta d_j^{\nu+1} = d_j^{\nu+1} - d_j^{\nu}, \label{1:update_k}\\
    & && \Delta m_j^{\nu+1} = m_j^{\nu+1} - m_j^{\nu}, \label{1:update_m}
\end{alignat}
\end{subequations}
with the objective to minimize the real part $(\sigma_\mathrm{max})$ of the rightmost (i.e., the most unstable) eigenvalue at each iteration step, and $[m^\mathsf{T},d^\mathsf{T}]^\mathsf{T}\in\R_{\geq0}^{2n_j}$ being the vector of decision variables. Constraint \eqref{1:sigma_iter} iteratively computes the real parts of all modes based on their previous values and the updates arising from incremental step changes in $m$ and $d$, while \eqref{1:sigma_lim} is needed for achieving the aforementioned objective. Inequalities \eqref{1:k_lim}-\eqref{1:m_lim} impose upper and lower bounds on total inertia and damping at each node, whereas \eqref{1:step_k}-\eqref{1:step_m} place limits on the permissible changes at each iteration; $\phi^d\in\R^{n_j}$ and $\phi^m\in\R^{n_j}$ represent the normalization of step size limits based on parameter sensitivities, which will be further elaborated in Section~\ref{subsec:soln_strategy}. Finally, \eqref{1:update_k}-\eqref{1:update_m} declare the updated decision variables for the next iteration step. The optimization ends once all modes become stable, i.e., the condition $\sigma_i<0,\forall i\in \N_{\leq n}$ is met. 

\subsubsection{Step 2: Improving Damping Ratio}

Once the system is small-signal stable, the next step is to make sure that the worst-case damping ratio is above a predefined threshold $\ushort{\zeta}\in\R_{>0}$. Apart from improving the damping ratio, limits are placed on RoCoF and frequency nadir to ensure an acceptable frequency response, which leads to the following problem:
\begin{subequations} \label{opti_damping}
\begin{alignat}{3}
    & \underset{\Phi_2, \eta}{\min} \quad && - c_\zeta \zeta_\mathrm{min} + c_f \big(\eta_{f_1}+\eta_{f_2}\big) + c_{\dot f} \big(\eta_{\dot{f}_1}+\eta_{\dot{f}_2}\big) \label{2:obj}\\
    & \;\mathrm{s.t.} \quad && \forall j \in \mathcal{J}, \forall i \in \N_{\leq n}, \nonumber \\
    & &&\eqref{1:k_lim}\text{-}\eqref{1:update_m}, \nonumber\\
    & && \zeta^{\nu+1}_i = \zeta^{\nu}_i + \sum_{j \in \mathcal{J}}{\frac{\partial \zeta^{\nu}_i}{\partial d_j}} \Delta d_j^{\nu+1} + \sum_{j \in \mathcal{J}}{\frac{\partial \zeta^{\nu}_i}{\partial m_j}} \Delta m_j^{\nu+1}, \label{2:zeta_update}\\
    & && \zeta_\mathrm{min} \leq \zeta^{\nu+1}_i, \label{2:zeta_min}\\
    & && D^{\nu+1} = \frac{\sum_{k\in\mathcal{K}}P_{g_k} d_k+\sum_{j\in\mathcal{J}}P_{g_j} d_j^{\nu+1}}{\sum_{n\in\mathcal{N}}P_{g_n}}, \label{2:D}\\
    & && M^{\nu+1} = \frac{\sum_{k\in\mathcal{K}}P_{g_k} m_k+\sum_{j\in\mathcal{J}}P_{g_j} m_j^{\nu+1}}{\sum_{n\in\mathcal{N}}P_{g_n}}, \label{2:M}\\
    & && \Delta D^{\nu+1} = D^{\nu+1} - D^{\nu}, \label{2:update_D}\\
    & && \Delta M^{\nu+1} = M^{\nu+1} - M^{\nu}, \label{2:update_M}\\
    & && \Delta f_\mathrm{max}^{\nu+1} = \Delta f_\mathrm{max}^{\nu} + \frac{\partial \Delta f_\mathrm{max}^\nu}{\partial D} \Delta D^{\nu +1} + \frac{\partial \Delta f_\mathrm{max}^\nu}{\partial M} \Delta M^{\nu +1}, \label{2:Nadir} \\
    & && -{\Delta \widebar{f}}_\mathrm{lim} - \eta_{f_1} \leq \Delta f_\mathrm{max}^{\nu+1} \leq {\Delta \widebar{f}}_\mathrm{lim} + \eta_{f_2}, \label{2:Nadir1} \\
    & && - \widebar{\dot f}_\mathrm{lim} - \eta_{\dot{f}_1} \leq f_0 \frac{\Delta P}{M^{\nu+1}} \leq \widebar{\dot f}_\mathrm{lim} +  \eta_{\dot{f}_2}, \label{2:RoCoF} \\
    & && \frac{M^{\nu+1}}{T} - F_g < D^{\nu+1}, \label{2:tm}\\
    & && \eta_{f_1},\eta_{f_2},\eta_{\dot{f}_1},\eta_{\dot{f}_2} \geq 0. \label{2:Slacks}
\end{alignat}
\end{subequations}
Here, $\Phi_2 \coloneqq [\Phi_1^\mathsf{T},\zeta^{\nu+1},M^{\nu+1},D^{\nu+1},\Delta M^{\nu+1},\Delta D^{\nu+1},\Delta f^{\nu+1}]^\mathsf{T}$ and $\eta \coloneqq [\eta_{f_1}, \eta_{f_2}, \eta_{\dot{f}_1}, \eta_{\dot{f}_2}]^\mathsf{T}$ are the vector of optimization and slack variables respectively, and $\zeta_\mathrm{min}$ represents the worst-case damping ratio of the system to be maximized, with $c_\zeta\in\R_{>0}$ being the corresponding cost factor. Equality \eqref{2:zeta_update} defines the new damping ratios of all modes based on their previous values and the corresponding updates of $m$ and $k$, whereas \eqref{2:zeta_min} ensures achieving the targeted objective. Constraints \eqref{1:k_lim}-\eqref{1:update_m} from the first step still apply, with additional expressions \eqref{2:D}-\eqref{2:Slacks} imposing limits on frequency metrics of interest discussed in Section~\ref{subsec:freq_metrics}. In particular, \eqref{2:D}-\eqref{2:update_M} define the total system inertia and damping as well as the incremental changes between iterations, \eqref{2:Nadir} describes the Taylor approximation of the nonlinear frequency nadir constraint, \eqref{2:RoCoF}-\eqref{2:Nadir1} provide upper and lower bounds on permissible RoCoF and frequency nadir magnitudes, and \eqref{2:tm} ensures that the time instance of frequency nadir is positive. Note that \eqref{2:RoCoF}-\eqref{2:Nadir1} are implemented as soft constraints, with slack variables defined by \eqref{2:Slacks} and included in \eqref{2:obj}, penalized by factors $c_f\in\R_{>0}$ and $c_{\dot{f}}\in\R_{>0}$. This optimization is completed when the the damping ratios of all modes reach a given threshold, i.e., $\zeta_\mathrm{min}\geq\ushort{\zeta}$, while also ensuring that the RoCoF and nadir conditions are satisfied.

\subsubsection{Step 3: Reducing Control Effort}

After achieving a satisfactory dynamic performance in terms of metrics defined in Section~\ref{subsec:dyn_prop}, the goal of the final optimization step is to reduce the total amount of inertia and damping in the system without compromising the aforementioned dynamic performance. This can be interpreted as a reduction of the VSM control effort through redistribution of virtual inertia and damping among different power converters, and formulated by
\begin{subequations} \label{opti_MandD}
\begin{alignat}{3}
    & \underset{\Phi_2, \eta}{\min} \quad && c_M M + c_D D\label{3:obj}\\
    & \;\mathrm{s.t.} \quad && \forall j \in \mathcal{J}, \forall i \in \N_{\leq n}, \nonumber \\
    & &&\eqref{1:k_lim}\text{-}\eqref{1:update_m}, \eqref{2:zeta_update}\text{-}\eqref{2:tm}, \nonumber\\
    & && \zeta_i^{\nu+1} \geq \ushort{\zeta}, \label{3:zeta_th} \\
    & && \eta_{f_1},\eta_{f_2},\eta_{\dot{f}_1},\eta_{\dot{f}_2} = 0,
\end{alignat}
\end{subequations}
where $c_M\in\R_{>0}$ and $c_D\in\R_{>0}$ reflect the ``price'' of control gains, \eqref{1:k_lim}-\eqref{1:update_m} and \eqref{2:zeta_update}-\eqref{2:tm} encompass the system constraints from previous steps, and \eqref{3:zeta_th} ensures that the damping ratios stay above the minimum permissible limit. In contrast to \eqref{opti_damping}, the frequency constraints in \eqref{opti_MandD} are implemented as hard constraints, i.e., slack variables $\eta_{f_1},\eta_{f_2},\eta_{\dot{f}_1},\eta_{\dot{f}_2}$ are set to zero, and normalization of step sizes is neglected, i.e., $\phi_j^d=\phi_j^m=1,\forall j\in\mathcal{J}$. The optimization ends when the incremental change in aggregate values of inertia and damping between $5$ consecutive iteration steps reaches a predefined lower bound $\epsilon=10^{-4}$. 

\subsection{Uniform Optimization Problem}
An alternative approach is to combine the objectives and constraints of all three steps into a single ``uniform'' formulation, which could potentially lead to a more effective allocation of virtual gains. Indeed, while employing different and independent objectives for each stage of the multi-step method effectively identifies the optimal parameters (i.e., the local optimum) for that particular step, it can also yield a suboptimal final solution compared to the uniform formulation. For instance, in the first and second step (i.e., \eqref{1:obj} and \eqref{2:obj} respectively) there is no cost associated with the amount of inertia and damping being placed, which could lead to unnecessarily high allocation of virtual inertia and damping. Despite the reduction of the control effort in the third step, the final distribution of virtual control gains can end up significantly different compared to the uniform approach.

The uniform optimization problem can be formulated as follows:
\begin{subequations} \label{uniform}
\begin{alignat}{3}
    & \underset{\Phi_2, \eta, \eta_\zeta}{\mathrm{min}} && c_\zeta \eta_\zeta + c_f \big(\eta_{f_1} + \eta_{f_2} \big) + c_{\dot f} \big(\eta_{\dot{f}_1}+\eta_{\dot{f}_2}\big) + c_M M + c_D D \label{uniform:obj}\\
    & \,\mathrm{s.t.} \quad && \forall j \in \mathcal{J}, \forall i \in \N_{\leq n}, \nonumber \\
    & &&\eqref{1:k_lim}\text{-}\eqref{1:update_m}, \eqref{2:zeta_update}\text{-}\eqref{2:Slacks}, \nonumber\\
    & && \zeta_i^{\nu+1} + \eta_\zeta \geq \ushort{\zeta}, \label{uniform:slack_zeta}\\
    & && \eta_\zeta \geq 0. \label{uniform:Slack}
\end{alignat}
\end{subequations}
The expressions \eqref{1:k_lim}-\eqref{1:update_m} and \eqref{2:zeta_update}-\eqref{2:Slacks} include previously defined constraints, whereas \eqref{uniform:slack_zeta} introduces a relaxation of the minimum damping ratio requirement, with $c_\zeta\in\R_{>0}$ reflecting the cost of violating the respective limit and $\eta_\zeta\in\R_{\geq0}$ being the new slack variable. Such formulation ensures that once the criteria for the minimum damping ratio, RoCoF and frequency nadir are met, the appropriate slack variables become zero and stop affecting the objective function. At this point, the only non-zero terms in \eqref{uniform:obj} pertain to the control effort (i.e., virtual inertia and damping), which corresponds to the final step of the multi-step formulation. 

\subsection{Sensitivity Computation and Solution Strategy} \label{subsec:soln_strategy}
As mentioned previously, the eigensensitivities are nonlinear and an iterative approach is used to compute the updates. The system is linearized around the current parameter values to obtain the state-space model from which the damping ratios and their sensitivities are derived. For the simplified model used in \cite{Borsche2015}, it is possible to obtain the analytical expressions for the specific sensitivities of the damping ratio. However, when studying a realistic system with significantly higher level of detail, obtaining the linearized system model and the aforementioned sensitivities at each iteration is not straightforward and can result in a high computational burden. We overcome this issue by employing the Symbolic Math Toolbox in \textsc{Matlab}~\cite{MATLAB} and deriving a symbolic state-space representation of the system. Note that this is done only once, prior to initialization of the sequential program. The eigenvalues and their respective sensitivities are then numerically computed at each iteration and used to determine the optimal updates of virtual control gains.

%% Without Flowchart 

Due to the inherent nature of the problem, it is possible that during the course of optimization the update steps of the decision variables (i.e., control parameters) are too large, which would in turn imply that the linearization equilibrium and the sensitivities used to compute the updates are no longer valid. To mitigate this problem, the following strategy is adopted. Once the optimizer finds a solution, it returns the new values of the parameters along with the values of the damping ratios. The updated values are then used to obtain the linearized model of the system as well as the eigenvalues and their corresponding damping ratios, which are subsequently compared with the damping ratios obtained earlier from the solution of the optimizer. If the difference between the two is larger than a prescribed threshold, it indicates that the change in the parameters is too large for the current linearization and the step size is reduced by half. This process is repeated until the difference goes below a given threshold or until the step size becomes too small.

Another improvement to the algorithm is introduced by continuously readjusting the upper and lower bounds on incremental step changes $(\Delta \ushort{m}_j,\Delta \widebar{m}_j,\Delta \ushort{d}_j,\Delta \widebar{d}_j)$ between iterations, based on the damping ratio sensitivities represented by $\phi_j^d\in\R$ and $\phi_j^m\in\R$:
\begin{equation}
    \phi^m_j = \frac{\frac{\partial \zeta_\mathrm{min}}{\partial m_j}}{\underset{j \in \mathcal{J}}{\mathrm{max}}\frac{\partial \zeta_\mathrm{min}}{\partial m_j}}\,, \quad 
    \phi^d_j = \frac{\frac{\partial \zeta_\mathrm{min}}{\partial d_j}}{\underset{j \in \mathcal{J}}{\mathrm{max}}\frac{\partial \zeta_\mathrm{min}}{\partial d_j}}.
\end{equation}
Such procedure assigns larger step limits to parameters with a greater impact on the damping ratio. In particular, this ensures that the virtual inertia and damping is added or removed only to converters with a substantial influence on dynamic performance of the system.

\vspace{-0.2cm}
\section{Results} \label{sec:res}

\begin{figure}[!t]
    \centering
    \scalebox{0.45}{\includegraphics[]{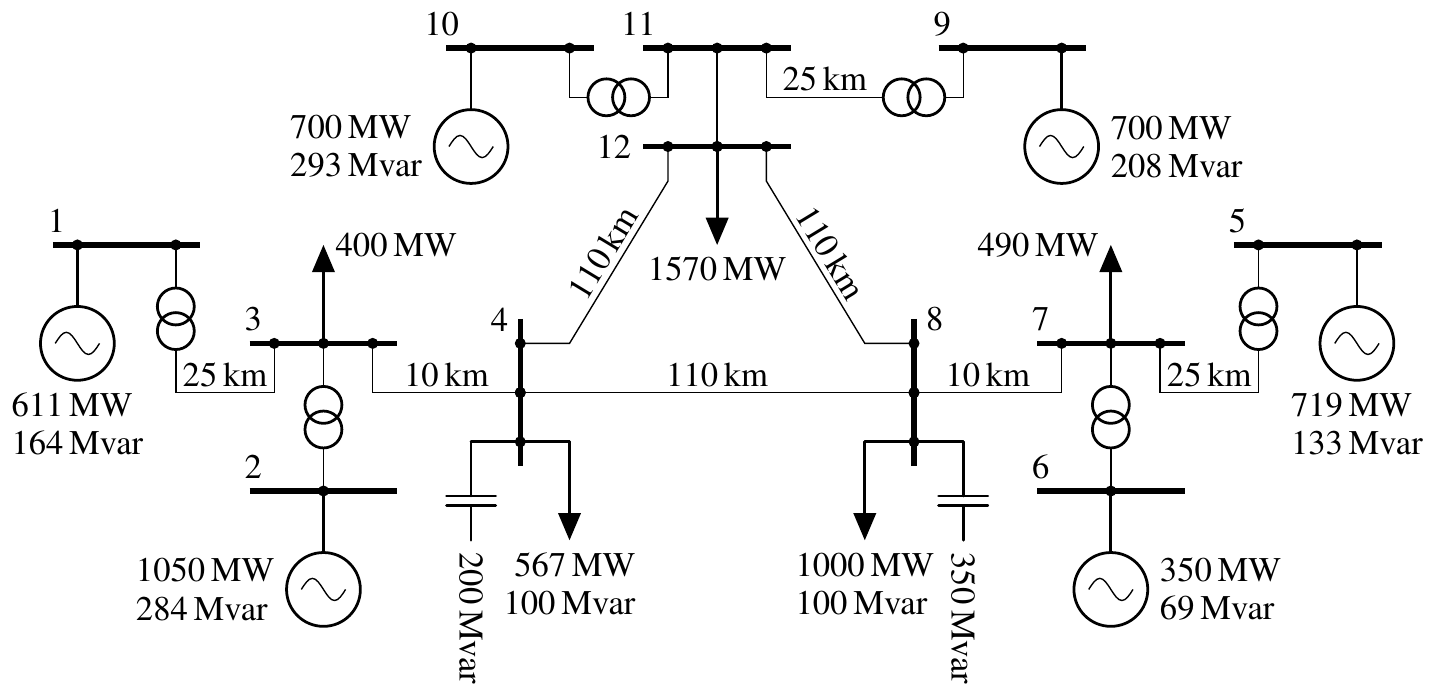}}
    \caption{Topology of the investigated three-area test system.}
    \label{fig:Kundur_3Area}
    \vspace{-0.35cm}
\end{figure}

We investigate the performance of the proposed optimization algorithms on a modified version of the well known Kundur's two-area system, with an addition of one more area forming a triangle depicted in Fig.~\ref{fig:Kundur_3Area}. The modified system comprises $6$ generators, with the network, generation and load parameters adapted from the original two-area system in \cite{Kundur1994}. The same test case has been previously used in other relevant studies on placement and effects of inertia and damping in low-inertia systems~\cite{Borsche2015,xiong2018stability,GrossIREP,ademola2018optimal}.
To emulate different system conditions, two test cases are considered: (i) a \textit{low-inertia} system with traditional SGs at nodes $1$ and $5$ and \textit{grid-following} converter-interfaced generation at the remaining four generation nodes; (ii) a \textit{no-inertia} system with a \SI{100}{\percent} inverter penetration, i.e., all six SGs replaced by VSC units. Nonetheless, only converters at nodes $1$ and $5$ are operating in \textit{grid-forming} mode, whereas the remaining VSCs are controlled as \textit{grid-following} units. As previously pointed out, we assume that only the inertia and damping constants of converters are controllable, with the initial VSM parameters taken as $m_j^{0}=\SI{0.5}{\second}$ and $d_j^{\,0}=2\,\mathrm{p.u.}$ The minimum permissible damping ratio is set at $\ushort{\zeta} = \SI{10}{\percent}$, whereas the thresholds enforced on maximum RoCoF and frequency deviation are $\widebar{\dot{f}}_\mathrm{max}=1\,\mathrm{Hz/s}$ and $\Delta\widebar{f}_\mathrm{max}=\SI{0.8}{\hertz}$. The system is subjected to a disturbance $\Delta P$, which is considered to be the worst-case power deficit caused by the loss of a single generator. For the purposes of this study it is assumed that the loss of generation occurs at node $1$. The upper and lower limits on incremental changes of inertia and damping at each iteration step are set to $\Delta \widebar{m}_j=\Delta \widebar{d}_j=0.5$ and $\Delta \ushort{m}_j=\Delta \ushort{d}_j=-0.5$, respectively.

Several different case studies are conducted. Firstly, a comparison between the virtual gain allocations obtained for the simplified model of a low-inertia system employed in \cite{BKP-SB-FD:17,xiong2018stability,GrossIREP,ademola2018optimal} and the detailed model from \cite{UrosStability} is presented. Secondly, the performance of the two proposed problem formulations on both test cases is investigated, together with the impact of frequency constraints on final distribution of virtual inertia and damping. Finally, the performance and convergence properties of the iterative algorithm from Section~\ref{subsec:soln_strategy} are discussed. 

\subsection{Simplistic vs Detailed Model} \label{subsec:res_model}

\begin{figure}[!b]
    \vspace{-0.35cm}
    \centering
    \scalebox{0.8}{\includegraphics[]{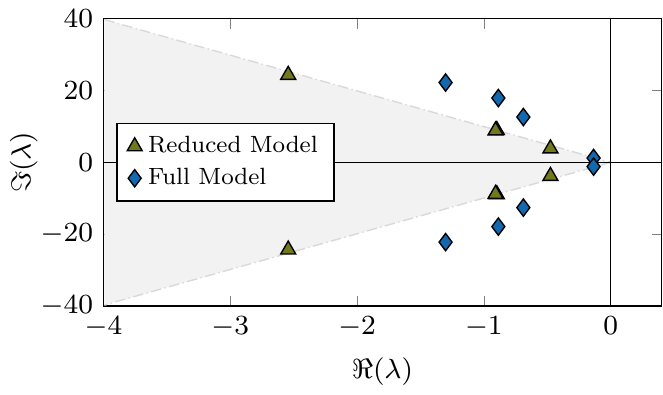}}
    \caption{Eigenvalue spectrum close to imaginary axis of the reduced- and full-order models for the identical inertia and damping allocation.}
    \label{fig:redFull_eig}   
\end{figure}

The goal of this case study is to demonstrate the need for a detailed system modeling when tackling the allocation of inertia and damping in a low-inertia system. First, the uniform optimization problem \eqref{uniform} is solved for the low-inertia test case described by the simplified model used in \cite{BKP-SB-FD:17,xiong2018stability,GrossIREP,ademola2018optimal}. Subsequently, the obtained allocation of inertia and damping is applied to the detailed model, with the corresponding eigenvalue spectrums of the most critical modes depicted in Fig.~\ref{fig:redFull_eig}. The shaded region indicates the root loci area where the damping ratio is above the predefined threshold of \SI{10}{\percent}. It is clear that the dynamic properties of a full-order model are not completely preserved in a simplified model, indicated by several modes having an unsatisfactory damping ratio. In other words, while results of the proposed algorithm meet the targeted objectives when applied to a reduced-order system, the achieved dynamic characteristics do not necessarily translate to a more realistic model, therefore diminishing the effectiveness and practicality of studies employing simplistic system representation used in \cite{Borsche2015}. 

%Moreover, similar discrepancies can be noticed between the optimal inertia and damping distribution for the two systems shown in Fig.~\ref{fig:redFull_bar}, both in terms of total allocation of virtual gains as well as the placement in the network. Note that the allocation corresponding to the full-order model in Fig.~\ref{fig:redFull_bar} does not correspond to the respective eigenvalue spectrum in Fig.~\ref{fig:redFull_eig}. More precisely, the distribution of inertia and damping in Fig.~\ref{fig:redFull_bar} is obtained by running the uniform optimization problem for each of the two model orders.

% \begin{figure}[!t]
%     \centering
%     \scalebox{0.9}{\includegraphics[]{figures/results/M_D_Reduced.pdf}}
%     \caption{Optimal distribution of inertia and damping through uniform formulation when applied to full- and reduced-order models.}
%     \label{fig:redFull_bar}  
% \end{figure}

% To be restored in the final paper

\subsection{Uniform vs Multi-Step Optimization} \label{subsec:sol_methods}

\begin{figure}[!t]
    \centering
    \scalebox{0.8}{\includegraphics[]{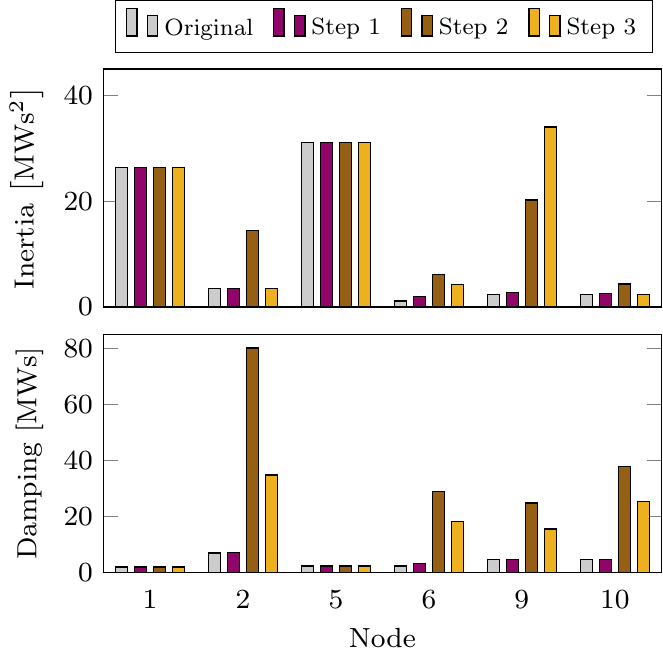}}
    \caption{Distribution of inertia and damping through multi-step optimization.}
    \label{fig:multi_step_bar}
    \vspace{-0.35cm}
\end{figure}

Here, we analyze the performance of the two optimization approaches presented in Section~\ref{sec:formulation}, applied to both low- and no-inertia test cases. To gain a better understanding of the sequential nature of the multi-step approach, we first study the outcome of each individual optimization step \eqref{opti_stability}-\eqref{opti_MandD}, with the respective distribution of inertia and damping illustrated in Fig.~\ref{fig:multi_step_bar} for the low-inertia case. 

The algorithm is capable of bringing the system to stable operation by adding small amounts of inertia and damping at nodes $6$ and $9$, justified by the fact that controller states associated with VSCs at those nodes have the highest participation in the unstable modes. However, the dynamic performance of the system does not meet prescribed requirements due to low damping ratios and high values of RoCoF and frequency nadir. Therefore, in the next step, a considerable amount of inertia and damping is placed at every node with converter-interfaced generation, which resolves the aforementioned issues. Finally, Step~$3$ reduces the total amount of virtual inertia and damping by readjusting the control gains of all four VSCs, while still meeting the necessary frequency and damping ratio criteria. The optimizer reduces inertia and damping at nodes that have low or even negative sensitivities and redistributes it to nodes with a larger impact on the system damping ratio. 

\begin{figure}[!t]
\centering
\begin{minipage}{1\columnwidth}
    \centering
    \hspace{0.05cm}
    \scalebox{0.8}{\includegraphics[]{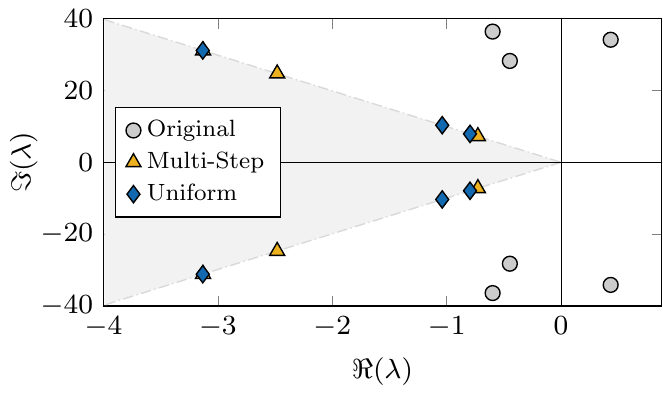}} 
    \vspace{0.5em}    
\end{minipage}
\begin{minipage}{1\columnwidth}
    \centering
    \scalebox{0.8}{\includegraphics[]{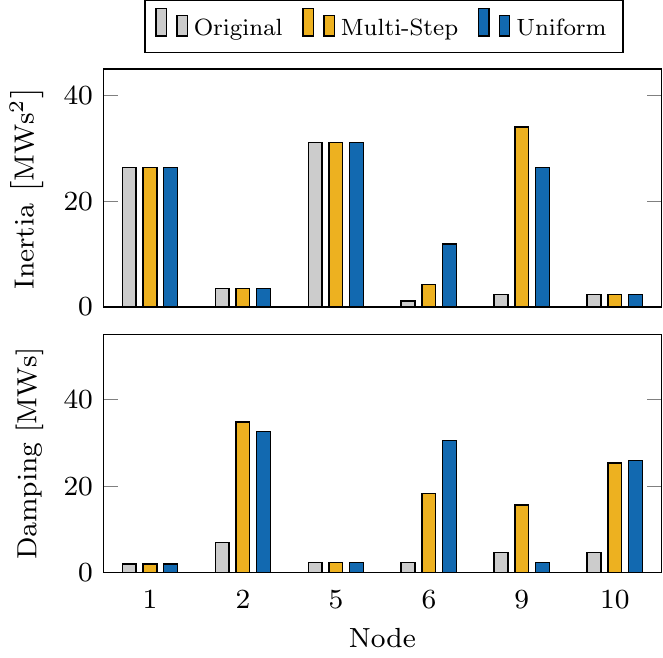}}   
\end{minipage} 
\caption{\label{fig:sol_SISIII}Comparison of uniform and multi-step formulation applied to a low-inertia system: critical eigenvalue spectrum (top) and allocation of inertia and damping (bottom).}
\vspace{-0.35cm}
\end{figure}

We can now compare the performance of the multi-step optimization to that of the uniform formulation. The eigenvalue spectrums in Fig.~\ref{fig:sol_SISIII} indicate that the system is unstable at the start of the optimization (corresponding to the \textit{original} allocation). Nevertheless, both formulations are capable of restoring stability and achieving satisfactory damping ratios and frequency response. Moreover, the total amount of virtual inertia and damping used by the two approaches remains the same, with the individual allocation differentiating between the two methods\footnote{The inertia and damping at nodes $1$ and $5$ remain intact due to the fact that only SGs are connected at these buses.}. This suggests that the frequency requirements, precisely RoCoF limit, act as binding constraints and impose a minimum aggregate inertia constant. Interestingly enough, in contrast to the multi-step approach and the addition of damping at all converter nodes, the uniform formulation increases virtual damping of only one VSC in each area, namely at nodes $2$, $6$ and $10$. In general, the final inertia placement is such that the aggregate inertia in different areas is approximately the same. This qualitatively matches the results obtained in \cite{Borsche2015} and \cite{GrossIREP}, where the distribution of inertia is similarly even across different areas. However, no such correlation can be made for the allocation of virtual damping. 

\begin{figure}[!t]
\centering
\begin{minipage}{1\columnwidth}
    \centering
    \hspace{0.175cm}
    \scalebox{0.8}{\includegraphics[]{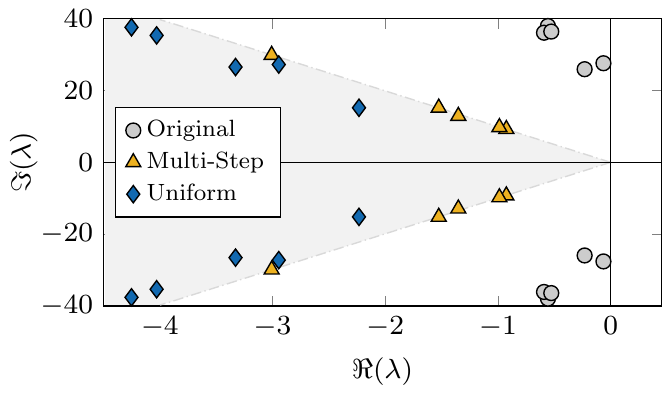}}
    \vspace{0.5em}    
\end{minipage}
\begin{minipage}{1\columnwidth}
    \centering
    \scalebox{0.8}{\includegraphics[]{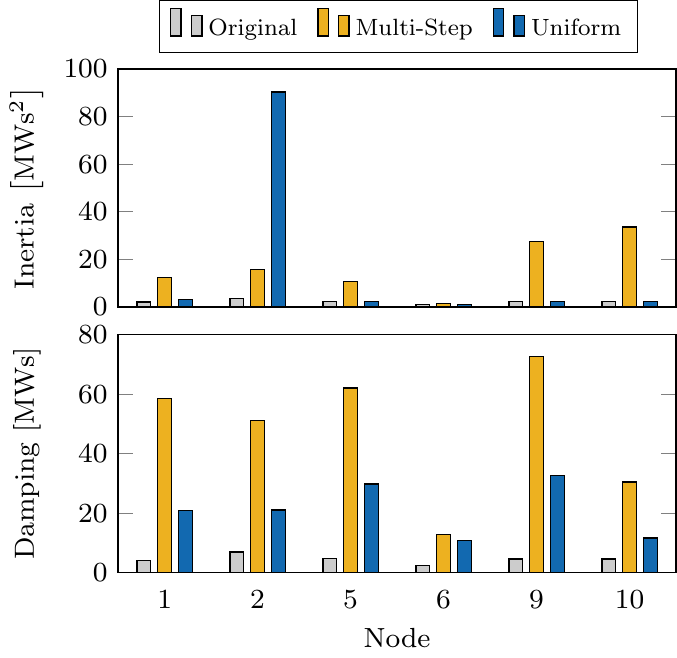}}   
\end{minipage} 
\caption{\label{fig:sol_IIIIII}Comparison of uniform and multi-step formulation applied to a no-inertia system: critical eigenvalue spectrum (top) and allocation of inertia and damping (bottom).}
\vspace{-0.35cm}
\end{figure}

Fig.~\ref{fig:sol_IIIIII} showcases the performance of the two algorithms under the no-inertia scenario, with the differences between the two allocations being more pronounced. Indeed, the multi-step approach yields a more even distribution of inertia but in turn employs more virtual damping. The underlying reason for such discrepancy lies in the problem formulation. In the multi-step approach, the cost for inertia and damping is only included in the final step, which resembles the most the optimization formulation of the uniform method. However, this also suggests that the starting point (and hence the computed sensitivities of interest) of the two algorithms will be different prior to the final step. In particular, due to no explicit cost for virtual control gains in the objective functions \eqref{1:obj} and \eqref{2:obj}, the solution of the first two sequences of the multi-step approach will reach a local optimum with a substantially higher installation of inertia and damping. On the other hand, the uniform formulation ensures that at each iteration only the minimum (i.e., necessary) amount of virtual gains is added to meet the prescribed system-level constraints. This might result in an uneven allocation of the parameters across the system, as shown in Fig.~\ref{fig:sol_IIIIII}, but it can easily be resolved by including an additional cost to promote the even distribution of inertia and damping. 

Another distinction in the outcomes of the two approaches is the placement of critical modes. In the multi-step approach, the critical eigenvalues are placed at the boundary of the shaded region, i.e., their damping ratios are close to the predefined threshold, whereas the uniform optimization results in critical modes being well within the shaded region. Due to the nature of the proposed problem formulation, the first few iterations of the multi-step approach focus on improving the frequency constraints, while the uniform method prioritizes the improvement of the damping ratio. Consequently, using the multi-step formulation the frequency criteria is met within fewer iterations compared to its uniform counterpart, and once the damping ratio is sufficiently high the optimization ends. In contrast, in the uniform approach the damping ratio criteria is satisfied first, with the optimization procedure continuing until the frequency limits (RoCoF in particular) are met. This results in higher final values of the damping ratios, as they are gradually increased over the course of remaining iterations. Nonetheless, such properties are not observed in the low-inertia test case (see Fig.~\ref{fig:sol_SISIII}) due to higher amount of inertia present in the system at the initialization stage, leading to lower RoCoF and frequency nadir values at the start of the optimization. 

Even though the uniform approach has a more generic problem formulation, there are applications for which the multi-step method would be more useful. For instance, the multi-step approach can identify the key parameters of interest for achieving different objectives. This information could also be used for potentially selecting other decision variables for different optimization steps. However, for the sake of simplicity, only the uniform approach will be studied in the remainder of the paper. 

\subsection{Impact of Frequency Constraints on Dynamic Performance} \label{subsec:freq_cons}

\begin{figure}[!b]
    \centering
    \vspace{-0.35cm}
    \centering
    \scalebox{0.8}{\includegraphics[]{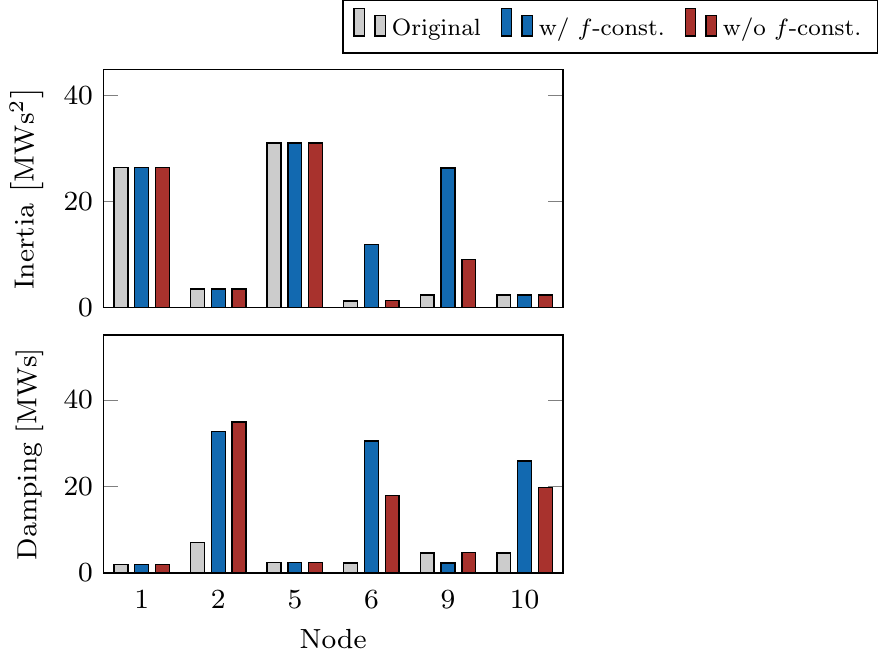}}\hspace{-2.5cm}
    \scalebox{0.8}{\includegraphics[]{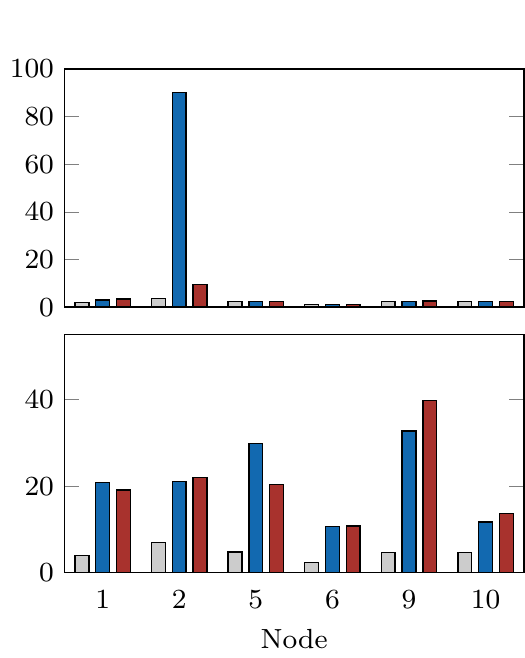}}
    \caption{\label{fig:freq_SISIII}Impact of frequency constraints on allocation in the low-inertia test case (left) and the no-inertia test case (right).}
\end{figure}

The impact of frequency-related constraints, namely the limits on maximum permissible RoCoF and frequency nadir, on the final solution is investigated by applying the uniform optimization problem to the same two test cases, with and without the inclusion of frequency constraints in \eqref{2:Nadir1}-\eqref{2:tm}. \begin{comment} Understandably, by removing the frequency constraints the critical modes in both scenarios are placed at the boundary of the region representing the set of acceptable damping ratios.\end{comment} 
It can be noticed in Fig.~\ref{fig:freq_SISIII} that the total amount of inertia and to some extent damping placed in the network is lower when the RoCoF and frequency nadir constraints are not considered and the reduction is more pronounced in the no-inertia test case. This is primarily due to the removal of the maximum RoCoF requirement, i.e., the minimum level of aggregate inertia needed in the network. In particular, the sensitivity of the minimum damping ratio to inertia and damping is in general of opposite sign, with the former being negative and latter being positive, which implies that the prescribed damping ratio criteria could be met by simultaneously reducing virtual inertia and damping gains at certain nodes. This indicates that in converter-dominated power systems the allocation of inertia is mostly influenced by the RoCoF constraint, whereas the damping distribution is more relevant for preserving the frequency nadir within thresholds and improving the critical damping ratio.

\subsection{Algorithm Performance and Convergence Properties} \label{subsec:perf}

The performance of the uniform method on a low-inertia test case, as well as the impact of imposing frequency constraints, can be studied in more detail by observing the appropriate metrics given in Table~\ref{tab:performance_metrics}. As previously discussed, the original system has low levels of aggregate inertia and damping and is small-signal unstable, indicated by the negative value of $\zeta_\mathrm{min}$. Additionally, the values of RoCoF and frequency nadir exceed the predefined thresholds and the system norms\footnote{Since the original system is unstable (resulting in $\| G \|_2=\| G \|_\infty=\infty$), the system norms for this case are computed after employing the first stage of the multi-step formulation, i.e., bringing the system to stability.} are particularly high, suggesting an unacceptable frequency response in case of a disturbance. Note that the $\mathcal{H}_\infty$ norm has been computed using the bisection method presented in \cite{boyd1989bisection}.

\begin{table}[!t]
\renewcommand{\arraystretch}{1}
\caption{Comparison of system performance metrics for the low-inertia test system and loss of SG at node $1$.}
\label{tab:performance_metrics}
\noindent
\centering
    \begin{minipage}{\linewidth} %Use the minipage environment to footnote tables
    \renewcommand\footnoterule{\vspace*{-5pt}} %to remove the horizontal rule above the table footnote
    \begin{center}
    \scalebox{1}{%
        \begin{tabular}{ c || c | c | c }
            \toprule
            \textbf{Metric} & \textbf{Original} & \textbf{w/ $\boldsymbol{f}$-const.} & \textbf{w/o $\boldsymbol{f}$-const.} \\ 
            \cline{1-4}
            Inertia $[\mathrm{MWs^2}]$ & $66.8$    & $101.5$    & $72.3$        \\
            Damping $[\mathrm{MWs}]$  & $23$    & $95.8$       & $62.9$        \\
            $\zeta_\mathrm{min}$       & $-0.01$  & $0.1$         & $0.1$          \\
            $\lvert \dot f_\mathrm{max} \rvert \, [\mathrm{Hz/s}]$   & $1.59$     & $1$           & $1.4$          \\
            $\lvert \Delta f_\mathrm{max} \rvert \, [\mathrm{Hz}]$   & $2.17$     & $0.58$        & $0.76$         \\
            $\mathcal{H}_2$ gain  & $18.8$ & $1.07$ & $1.87$         \\
            $\mathcal{H}_\infty$ gain   & $3.71$     & $0.76$        & $1.11$        \\
            \arrayrulecolor{black}\bottomrule
        \end{tabular}
    }
        \end{center}
    \end{minipage}
\vspace{-0.35cm}
\end{table}

By solving the optimization problem with frequency constraints, the total amount of inertia and damping in the system increases by \SI{52}{\percent} and \SI{316}{\percent}, respectively, which results in the worst-case damping ratio reaching the exact predefined threshold of \SI{10}{\percent}. Moreover, the values of RoCoF and nadir are now within their acceptable limits, with RoCoF being at the prescribed boundary. As a consequence, the $\mathcal{H}_2$ and $\mathcal{H}_\infty$ gains are substantially lower compared to their initial values. On the other hand, when the frequency constraints are not considered, the system inertia and damping increase by \SI{8}{\percent} and \SI{173}{\percent} respectively. While the total amount of employed virtual gains is lower in this case, the maximum RoCoF is unacceptably high and the $\mathcal{H}_2$ and $\mathcal{H}_\infty$ norms have larger values, thus suggesting a necessary trade-off between the control effort and dynamic performance.

\begin{figure}[!t]
    \centering
    \hspace{-0.25cm}
    \scalebox{0.0485}{\includegraphics[]{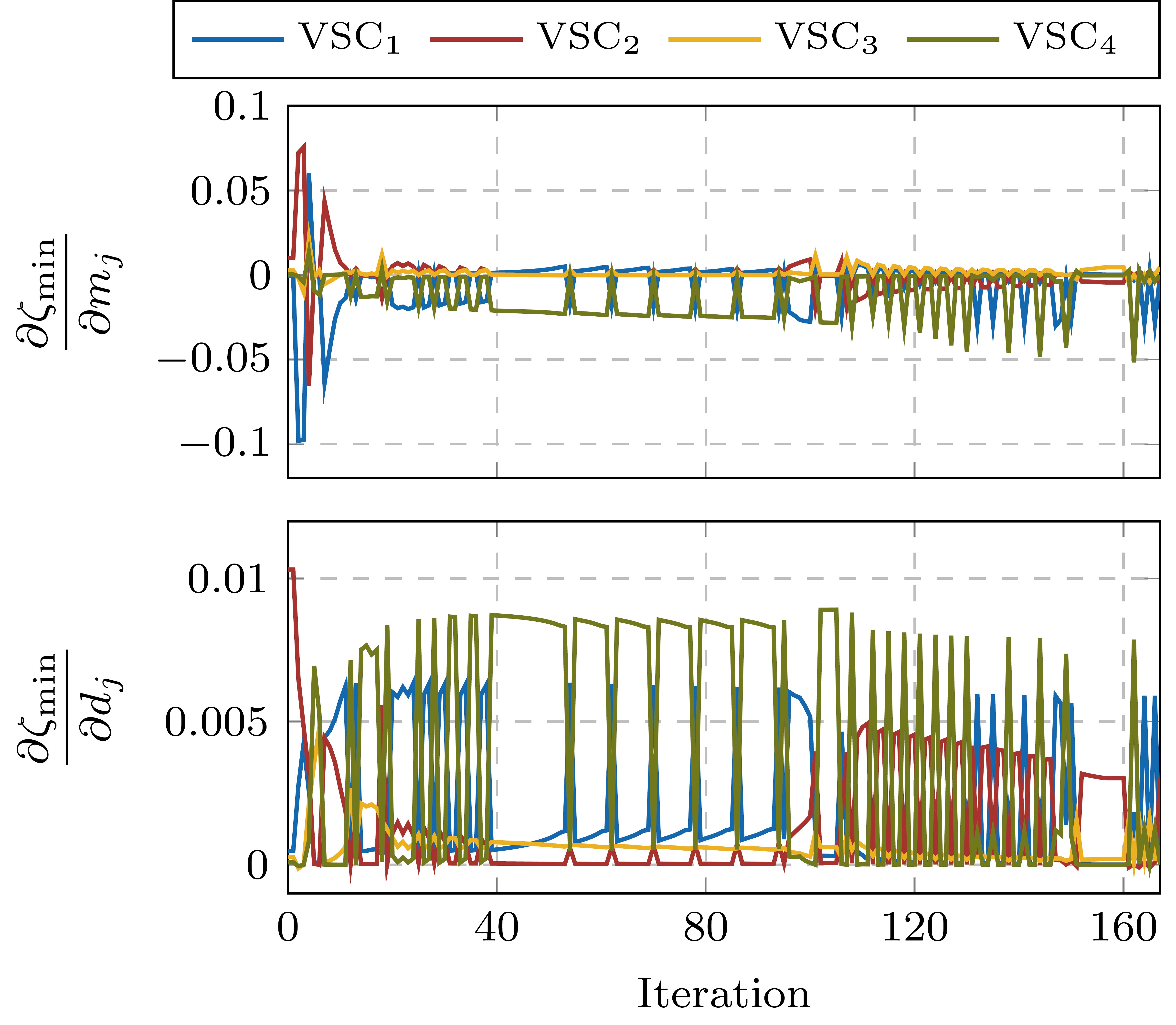}}
    \caption{Sensitivity evolution of the worst-case damping ratio to virtual inertia and damping parameters of individual converters over iterations.}
    \label{fig:sensitivity_plot}
    \vspace{-0.35cm}
\end{figure}

Finally, we study the numerical characteristics and convergence of the proposed algorithm. The sensitivities of the worst-case damping ratio with respect to virtual inertia and damping constants of converters are highly nonlinear, illustrated by the evolution of respective sensitivities over iterations in Fig.~\ref{fig:sensitivity_plot}. The oscillatory behavior justifies the need for adaptive step sizing of parameter updates and the techniques mentioned in Section~\ref{subsec:soln_strategy}. Another key inference is that the sensitivity of the damping ratio to inertia is in general negative and smaller in magnitude compared to the positive sensitivity to damping, which further supports the claim that damping is a more relevant control gain of the two for improving the damping ratios in the system. 

% \begin{figure}[!t]
%     \centering
%     \scalebox{0.8}{\includegraphics[]{figures/results/damping_H2.pdf}}
%     \caption{Evolution of the worst-case damping ratio and system $\mathcal{H}_2$ norm over iterations.}
%     \label{fig:damp_h2_movement}
%     \vspace{-0.35cm}
% \end{figure}

% \begin{figure}[!t]
%     \centering
%     \hspace{-0.25cm}
%     \scalebox{0.8}{\includegraphics[]{figures/results/sensitivity_plot.pdf}}
%     \caption{Sensitivity evolution of the worst-case damping ratio to virtual inertia and damping parameters of individual converters over iterations.}
%     \label{fig:sensitivity_plot}
%     \vspace{-0.35cm}
% \end{figure}

% \begin{figure}[!b]
%     \centering
%     \vspace{-0.35cm}
%     \hspace{0.25cm}
%     \scalebox{0.8}{\includegraphics[]{figures/results/damping_H2_new.pdf}}
%     \caption{Evolution of the worst-case damping ratio and $\mathcal{H}_2$ norm over iterations.}
%     \label{fig:damp_h2_movement}
% \end{figure}

\begin{figure}[!b]
    \centering
    \vspace{-0.35cm}
    \hspace{0.25cm}
    \scalebox{0.0485}{\includegraphics[]{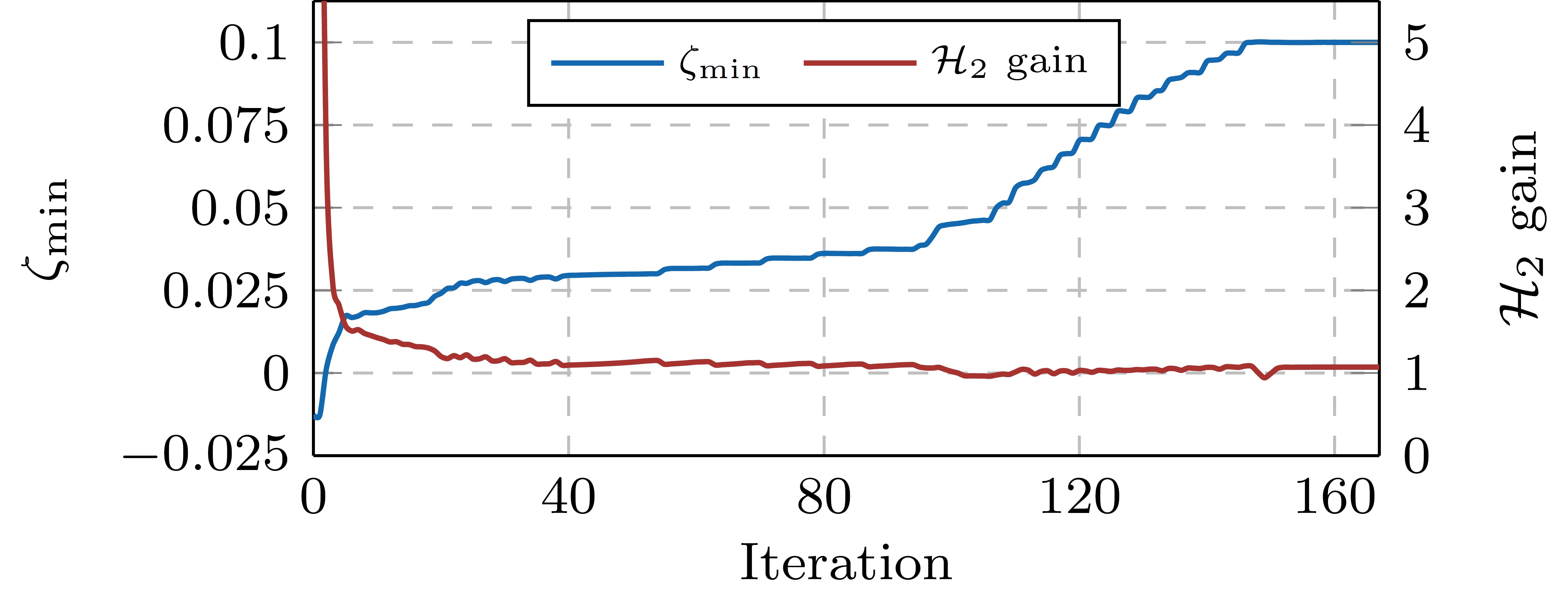}}
    \caption{Evolution of the worst-case damping ratio and $\mathcal{H}_2$ norm over iterations.}
    \label{fig:damp_h2_movement}
\end{figure}

The progressive iterative improvement of the worst-case damping ratio and the value of the $\mathcal{H}_2$ gain during the course of optimization are presented in Fig.~\ref{fig:damp_h2_movement}. An important observation is that the rate-of-change of the worst-case damping ratio varies considerably as a direct consequence of the oscillatory nature of sensitivities from Fig.~\ref{fig:sensitivity_plot}. Moreover, the increase in the damping ratio and the reduction of RoCoF and frequency nadir lead to a more desirable frequency response, which in turn contributes to the improvement of system norms. This validates the claims pertaining to similarities between improving the damping ratio and system norms raised in Section~\ref{sec:problem}, since maximizing the damping ratio while simultaneously limiting the frequency response metrics achieves a similar target as directly minimizing the $\mathcal{H}_2$ and $\mathcal{H}_\infty$ norms in \cite{GrossIREP}. However, the proposed approach imposes lower computational requirements and can be directly applied to realistic low-inertia systems. Further inspection of Fig.~\ref{fig:damp_h2_movement} reveals that during the last $15$ iterations the worst-case damping ratio and the $\mathcal{H}_2$ norm are fairly constant. This segment indicates the process of inertia and damping redistribution, i.e., the minimization of the employed control effort. 

\section{Conclusion} \label{sec:conclusion}

In this paper, we improve the dynamic performance of an inverter-dominated power system by optimally allocating virtual inertia and damping across the network. We tackle this problem by formulating an optimization problem based primarily on the sensitivities of damping ratios to inertia and damping constants of individual generators. Moreover, we consider several additional performance metrics such as the frequency nadir, RoCoF and small-signal stability, and incorporate them as explicit constraints into two conceptually different iterative problem formulations. Furthermore, improvements in terms of accounting for the multi-objective nature of the problem, computational efficiency and adaptive step-size adjustments have also been made. 

The results indicate that the simplified system models, commonly used in the literature, do not accurately capture the dynamics of a power system with both conventional and converter-interfaced generation. In other words, a detailed representation of low-inertia grids is needed when dealing with inertia and damping allocation problems, which poses issues for the existing methods based on minimizing system norms. While conceptually different, both proposed formulations provide meaningful results and insightful observations in terms of the overall impact of different VSM control gains and frequency-related constraints on system dynamics. However, due to its more generic formulation and multifaceted objective function, we conclude that the uniform approach is the preferred method of the two. Moreover, we show that by improving the worst-case damping ratio and constraining the frequency metrics of interest, the algorithm also achieves a significant reduction in $\mathcal{H}_2$ and $\mathcal{H}_\infty$ norms, therefore combining the objectives and targets of various studies in the literature within a single optimization problem.

% References section
\bibliographystyle{IEEEtran}
\bibliography{thesis}
% That's all folks
\end{document}